\documentclass[a4paper,12pt]{article}
\usepackage{amsmath,latexsym,amsfonts,amssymb,amsthm,mathrsfs,longtable,lscape,dsfont,fancybox,yfonts}
\usepackage{amsthm}
\usepackage{amssymb}
\usepackage{amscd}
\usepackage{amsfonts}
\usepackage{amsbsy}
\usepackage{fancyhdr}
\usepackage{graphicx}
\usepackage[all]{xy}
\usepackage{epsfig}
\newtheorem {teo} {Theorem} [section]

\newtheorem {prop} [teo]{Proposition}

\title{Generalization of the Hill's problem as an application for the Trojan asteroids of the solar system}
\author{Jaime Burgos--Garc\'\i a \thanks{Yeshiva University, Beren Campus. 245 Lexington avenue, 10016. New York City. N.Y.
e--mail: jbg84@xanum.uam.mx} \and Marian
Gidea\thanks{Yeshiva University. Preprint.}}

\date{}

\begin{document}

\maketitle

\begin{abstract}
The restricted four-body problem studies the dynamics of a massless particle under the gravitational force produced by three masses (primaries) in an equilateral configuration. One primary, say $m_{3}$, is considered too small compared with the other ones. In a similar way as in the classical Hill's problem, we study the limit case $m_{3}\rightarrow0$ in the Hamiltonian of the R4BP. In this paper we prove that such limit exists and the resulting limit problem produces a new Hamiltonian that inherits some basic features of the restricted three and four body problems. We analyze some dynamical aspects of this new system that can be considered as a generalization of the Hill's problem.
\end{abstract}



\section{Introduction}
Few bodies problems have been studied for long time in celestial
mechanics,  either as simplified models of more complex planetary
systems or as benchmark models where new mathematical theories can
be tested. The three--body problem has been a source of inspiration
and study in Celestial Mechanics since Newton and Euler, in particular the restricted three body problem (R3BP) has
demonstrated to be a good model of several systems in our solar
system such as the Sun--Jupiter--Asteroid system, and with less
accuracy the Sun--Earth--Moon system, in these systems the R3BP was used to know preliminary orbits in some space missions. In analogy with the R3BP, in this paper we study a restricted problem of four bodies consisting
of three primaries moving in circular orbits keeping an equilateral
triangle configuration and a massless particle moving under the
gravitational attraction of the primaries.  It is known that in our solar system we can find such configurations, the so called Trojan asteroids of Jupiter, Mars and Neptune form approximately an equilateral configuration with their respective planet and the sun, Saturn--Tethys--Telesto, Saturn--Tethys--Calypso or Saturn--Dione--Helen are good examples of such configuration. Several authors \cite{PapaIII} ,\cite{Cecc}, have considered the restricted four body body problem to model the dynamics of a spacecraft in the Sun-Jupiter-Asteroid-spacecraft system. \\ \\G. W. Hill developed his famous lunar theory \cite{Hill} as an alternative approach for the study of the motion of the moon. As a first approximation, this approach consider a Kepler problem (Earth-Moon) with a gravitational perturbation produced by a far away massive body (Sun), some orbital elements such as the eccentricities of the orbits of the Moon and the earth and the inclination of the Moon are supposed to be zero. Previously to the Hill's work, the approach to study the dynamics of the moon consisted on considering two Kepler problems, one for the motion of the Earth and the Moon around their center of mass and other for the motion of the sun and such center of mass. However, this approach had several difficulties because of the solutions were given in terms of formal power series of orbitals elements, the principal inconvenience was due to the poor convergency of these series in terms of the ratios of the mean motions of the Earth and the Moon, the so called critical parameter. The success of the Hill's approach was given by using his model to obtain a periodic orbit of the trajectory of the Moon and then he included the orbital elements to correct it, in such a way, he avoided the computation of expansions in terms of the critical critical parameter. In a four body problem context, the smallness of one primary creates complicated equations of motion where an analytical study is extremely difficult to make and even there are technical inconveniences in the accuracy of numerical simulations. In the next sections we develop a model as a first approximation of the dynamics of a masses particle in a Sun-Planet-Asteroid system, as possible applications of this model we can consider the massless body like a spacecraft or a small satellite like the moon of the Trojan asteroid 624 Hektor \cite{Marchis}. In future works we may include relevant effects produced by inclinations and librations of the asteroids or perturbations due to other bodies for example.

\section{The restricted four body problem}
Consider three point masses, called $\textit{primaries}$, moving in
circular periodic orbits around their center of mass under their
mutual Newtonian gravitational attraction, forming an equilateral
triangle configuration. A fourth massless particle is moving under the gravitational attraction of the primaries, this problem is known as the equilateral restricted four body problem or simply as the restricted four body problem (R4BP). The equations
of motion in the usual dimensionless coordinates of the massless particle referred to a synodic frame of reference, where the primaries remain fixed, are:

\begin{figure}[!hbp]
\centering
\includegraphics[width=2.5in]{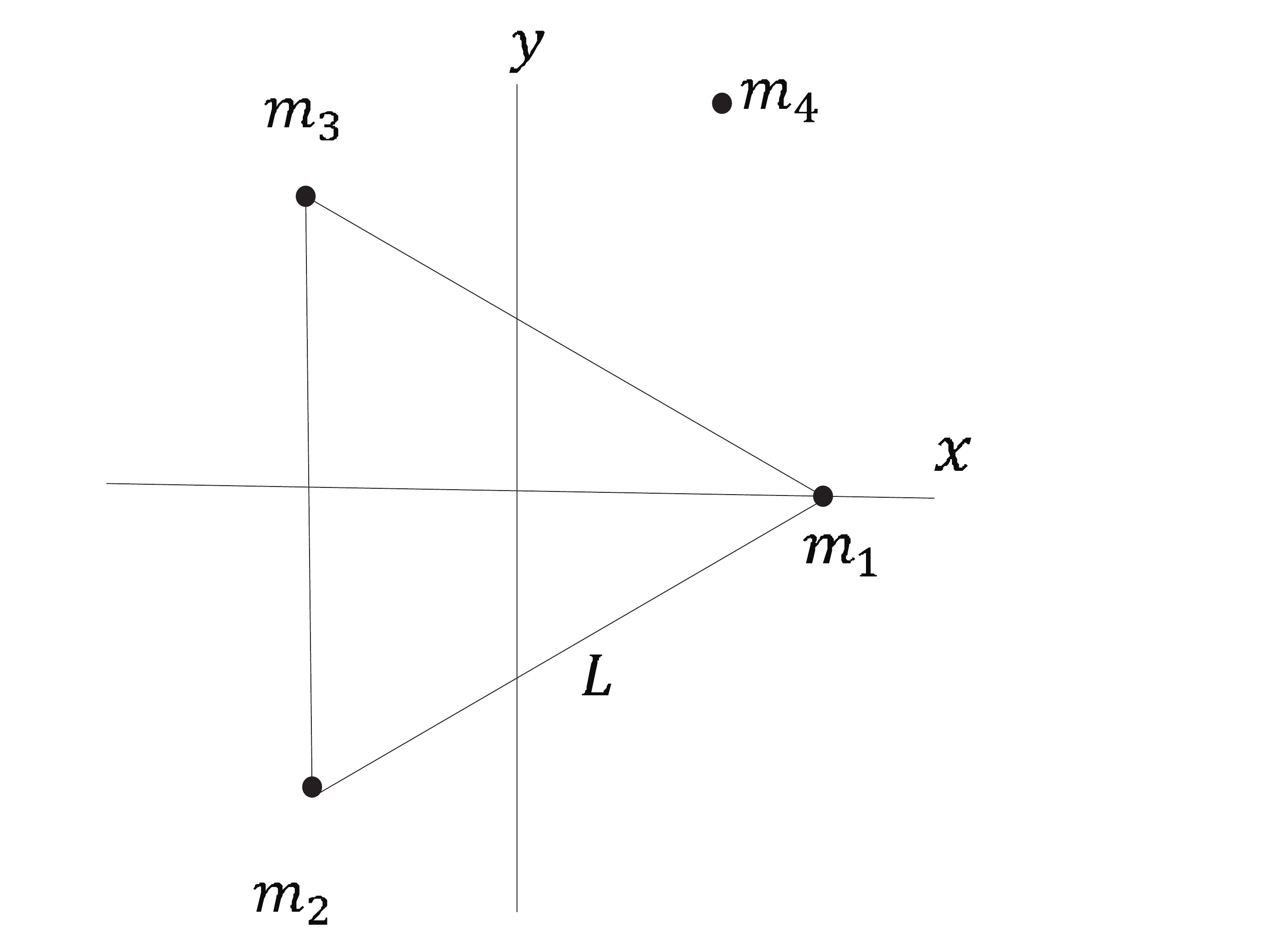}
\caption{The restricted four-body problem in a synodic system for the two equal masses case.\label{triangle}}
\end{figure}

\begin{equation}
\label{ecuacionesfinales}
\ddot{x}-2\dot{y}=\Omega_{x},
\end{equation}
$$\ddot{y}+2\dot{x}=\Omega_{y},$$ $$\ddot{z}=\Omega_{z},$$
where $$\Omega(x,y,z,m_{1},m_{2},m_{3})=\frac{1}{2}(x^{2}+y^{2})+\sum_{i=1}^{3}\frac{m_{i}}{r_{i}},$$
and $r_{i}=\sqrt{(x-x_{i})^{2}+(y-y_{i})^{2}+z^2}$, for $i=1,2,3$. The general expressions of the coordinates of the primaries in terms of the masses of the three primaries are given by
$$x_{1}=\frac{-\vert K\vert\sqrt{m_{2}^{2}+m_{2}m_{3}+m_{3}^{2}}}{K},$$ $$y_{1}=0,$$ $$x_{2}=\frac{\vert K\vert[(m_{2}-m_{3})m_{3}+m_{1}(2m_{2}+m_{3})]}{2K\sqrt{m_{2}^{2}+m_{2}m_{3}+m_{3}^{2}}},$$ 
\begin{equation}
y_{2}=\frac{-\sqrt{3}m_{3}}{2m_{2}^{3/2}}\sqrt{\frac{m_{2}^{3}}{m_{2}^{2}+m_{2}m_{3}+m_{3}^{2}}},\label{coordinatesprimaries}
\end{equation}
$$x_{3}=\frac{\vert K\vert}{2\sqrt{m_{2}^{2}+m_{2}m_{3}+m_{3}^{2}}},$$ $$y_{3}=\frac{\sqrt{3}}{2\sqrt{m_{2}}}\sqrt{\frac{m_{2}^{3}}{m_{2}^{2}+m_{2}m_{3}+m_{3}^{2}}}.$$

Where $K=m_{2}(m_{3}-m_{2})+m_{1}(m_{2}+2m_{3})$ and the three masses satisfy the relation $m_{1}+m_{2}+m_{3}=1$. It can be proved that the equations of motion have a first integral $$\dot{x}^2+\dot{y}^2+\dot{z}^2=2\Omega-C,$$ where $C$ is a constant. It is worth noting that when we make $m_{3}=0$ and $m_{2}:=\mu$ we recover the coordinates of the restricted three body problem (R3BP) $(x_{1},y_{1})=(-\mu,0)$, $(x_{2},y_{2})=(1-\mu,0)$ and $(x_{3},y_{3})=(1/2-\mu,\sqrt{3}/2)$, now the "phantom" mass $m_{3}$ is located in the so called equilibrium point $L_{4}$ of the R3BP. In the following, it will be necessary to consider the Hamiltonian of the system
\begin{equation}
H=\frac{1}{2}(p_{x}^{2}+p_{y}^{2}+p_{z}^{2})+yp_{x}-xp_{y}-\frac{m_{1}}{r_{1}}-\frac{m_{2}}{r_{2}}-\frac{m_{3}}{r_{3}}.\label{originalhamiltonian}
\end{equation}

\section{The limit case and equations of motion}
In this section we will discuss the how to compute the limit when $m_{3}\rightarrow0$ for the R4BP. We use a similar procedure as shown in \cite{MeyerHDS} by considering a symplectic scaling of the Hamiltonian and expansions in Taylor series in a neighborhood of the small mass $m_{3}$. The resulting Hamiltonian will be a three degrees of freedom system depending on a parameter $\mu$ which is the mass of the primary $m_{2}$.

\begin{teo} The limit $m_{3}\rightarrow0$ of the Hamiltonian (\ref{originalhamiltonian}) restricted to a neighborhood of $m_{3}$ exists and gives rise to a new Hamiltonian
\begin{equation}
H=\frac{1}{2}(p^{2}_{x}+p^{2}_{y}+p_{z}^{2})+yp_{x}-xp_{y}+\frac{1}{8}x^2-\frac{3\sqrt{3}}{4}(1-2\mu)xy-\frac{5}{8}y^2+\frac{1}{2}z^2\label{hillhamiltonian}
\end{equation}
$$-\frac{1}{\sqrt{x^2+y^2+z^2}},$$ where $m_{1}=1-\mu$ and $m_{2}:=\mu$.
\end{teo}

\textit{Proof}. We consider the Hamiltonian of the restricted four body problem (R4BP) in the center of mass coordinates
$$
H=\frac{1}{2}(p_{x}^{2}+p_{y}^{2}+p_{z}^{2})+yp_{x}-xp_{y}-\frac{m_{1}}{r_{1}}-\frac{m_{2}}{r_{2}}-\frac{m_{3}}{r_{3}},
$$ where $r_{i}^{2}=(x-x_{i})^2+(y-y_{i})^2+z^2$ and $(x_{i},y_{i})$ denotes the position of the primary $m_{i}$ for $i=1,2,3.$. We make the change of coordinates $x\rightarrow x+x_{3}$, $y\rightarrow y+y_{3}$, $z\rightarrow z$, $p_{x}\rightarrow p_{x}-y_{3}$, $p_{y}\rightarrow p_{y}+x_{3}$, $p_{z}\rightarrow p_{z}$, therefore in these new coordinates the Hamiltonian (\ref{originalhamiltonian}) becomes

\begin{equation}
H=\frac{1}{2}(p_{x}^{2}+p_{y}^{2}+p_{z}^{2})+yp_{x}-xp_{y}-(x_{3}x+y_{3}y)-\frac{m_{1}}{\bar{r}_{1}}-\frac{m_{2}}{\bar{r}_{2}}-\frac{m_{3}}{\bar{r}_{3}},\label{traslatedhamiltonian}
\end{equation}

where now we have $\bar{r}_{i}^{2}=(x+x_{3}-x_{i})^2+(y+y_{3}-y_{i})^2+z^2:=(x+\bar{x}_{i})^2+(y+\bar{y}_{i})^2+z^2$  for $i=1,2,3.$. We expand the terms $\frac{1}{\bar{r}_{1}}$ and $\frac{1}{\bar{r}_{2}}$ in Taylor series around the new origin of coordinates, if we ignore the constant terms we obtain the following expressions $$f^{1}:=\frac{1}{\bar{r}_{1}}=\sum_{k\ge1}P_{k}^{1}(x,y,z),$$ $$f^{2}:=\frac{1}{\bar{r}_{2}}=\sum_{k\ge1}P_{k}^{2}(x,y,z),$$ where $P_{k}^{j}(x,y,z)$ is a homogenous polynomial of degree $k$ for $j=1,2.$ In order to take the limit as $m_{3}\rightarrow0$, we perform the following symplectic scaling $x\rightarrow m_{3}^{1/3}x$, $y\rightarrow m_{3}^{1/3}y$, $z\rightarrow m_{3}^{1/3}z$, $p_{x}\rightarrow m_{3}^{1/3}p_{x}$, $p_{y}\rightarrow m_{3}^{1/3}p_{y}$ $p_{z}\rightarrow m_{3}^{1/3}p_{z}$ with multiplier $m_{3}^{-2/3}$, therefore 

\begin{equation}
H=\frac{1}{2}(p_{x}^{2}+p_{y}^{2}+p_{z}^{2})+yp_{x}-xp_{y}-\frac{1}{\bar{r}_{3}}-m_{3}^{-1/3}(x_{3}x+y_{3}y+P_{1}^{1}+P_{1}^{2})- \label{scaled}
\end{equation}

$$\sum_{k\ge2}m_{3}^{\frac{k-2}{3}}m_{1}P_{k}^{1}(x,y,z)-\sum_{k\ge2}m_{3}^{\frac{k-2}{3}}m_{2}P_{k}^{2}(x,y,z).$$ A straightforward computation shows $$P_{1}^{1}=m_{1}(\frac{x_{1}-x_{3}}{\bar{r}_{1}^{3}}x+\frac{y_{1}-y_{3}}{\bar{r}_{1}^{3}}y),$$ $$P_{1}^{2}=m_{2}(\frac{x_{2}-x_{3}}{\bar{r}_{2}^{3}}x+\frac{y_{2}-y_{3}}{\bar{r}_{2}^{3}}y),$$ where $\bar{r_{i}}=\sqrt{(x_{3}-x_{i})^2+(y_{3}-y_{i})^2}.$ It is important to note that the the first partial derivative is given by $$f^{i}_{z}=-\frac{z}{\bar{r}^{3}_{i}},$$ for $i=1,2.$ Therefore we obtain\begin{eqnarray*}
f^{i}_{z}(0,0,0)=f^{i}_{xz}(0,0,0)=f^{i}_{yz}(0,0,0)=0, \label{secondnullderivatives}
\end{eqnarray*}
and $$f^{i}_{zz}(0,0,0)=-1.$$ Now if we recall that the three masses are in equilateral configuration and we use the relation $m_{1}=1-m_{2}-m_{3}$ we obtain $$m_{3}^{-1/3}(x_{3}x+y_{3}y+P_{1}^{1}+P_{1}^{2})=m_{3}^{-1/3}[x_{1}+m_{2}(x_{2}-x_{1})-m_{3}(x_{1}-x_{3})]x$$ $$-m_{3}^{-1/3}[y_{1}+m_{2}(y_{2}-y_{1})-m_{3}(y_{1}-y_{3})]y,$$ in terms of the coordinates of the primaries (\ref{coordinatesprimaries}), we can write $$m_{3}^{-1/3}[y_{1}+m_{2}(y_{2}-y_{1})-m_{3}(y_{1}-y_{3})]=-m_{3}^{2/3}m_{2}s_{1}(m_{1},m_{2},m_{3})+m_{3}^{2/3}y_{3},$$ where $$s_{1}(m_{1},m_{2},m_{3})=\sqrt{\frac{3m_{2}^{3}}{4m_{2}^{3}(m_{2}^{2}+m_{2}m_{3}+m_{3}^{2})}},$$ and $$y_{3}=\frac{\sqrt{3}}{2\sqrt{m_{2}}}\sqrt{\frac{m_{2}^{3}}{m_{2}^{2}+m_{2}m_{3}+m_{3}^{2}}}.$$
A similar computation shows that the coefficient $m_{3}^{-1/3}[x_{1}+m_{2}(x_{2}-x_{1})-m_{3}(x_{1}-x_{3})]$ can be written in terms of a positive power of $m_{3}$. Therefore, the Hamiltonian (\ref{scaled}) looks like
\begin{equation}
H=\frac{1}{2}(p_{x}^{2}+p_{y}^{2}+p_{z}^{2})+yp_{x}-xp_{y}-\frac{1}{r}-m_{1}P_{2}^{1}-m_{2}P_{2}^{2}+\mathcal{O}(m_{3}^{1/3}).\label{unlimitedhamiltonian}
\end{equation}
We have defined $r=\bar{r}_{3}$. Now, we are allowed to take the limit $m_{3}\rightarrow0$ in the expression (\ref{unlimitedhamiltonian}), the final expression looks like
\begin{equation}
H=\frac{1}{2}(p^{2}_{x}+p^{2}_{y}+p_{z}^{2})+yp_{x}-xp_{y}+\frac{1}{8}x^2-\frac{3\sqrt{3}}{4}(1-2\mu)xy-\frac{5}{8}y^2+\frac{1}{2}z^2-K(x,y,z),\label{hillhamiltonian}
\end{equation}
where $K(x,y,z)=\frac{1}{\sqrt{x^2+y^2+z^2}}$, $m_{2}:=\mu$ and $m_{1}=1-\mu$. \qed

\textbf{Remarks}
\begin{itemize}
\item The expression 
\begin{equation}
Q=\frac{1}{2}(p^{2}_{x}+p^{2}_{y}+p_{z}^{2})+yp_{x}-xp_{y}+\frac{1}{8}x^2-\frac{3\sqrt{3}}{4}(1-2\mu)xy-\frac{5}{8}y^2+\frac{1}{2}z^2,
\end{equation}
is the quadratic part of the Hamiltonian of the restricted three body problem centered in the so called equilibrium point $L_{4}$.
\item Because of the above remark, we can consider the mass parameter in the range $\mu\in[0,1/2]$, the case where $\mu=1/2$ corresponds to the equal massive bodies case.
\item It will proved in the next section that this system have 4 equilibrium points in a neighborhood of $m_{3}$ and such equilibrium points will posses the same stability properties as in the full R4BP when $m_{3}$ is small but non zero.
\end{itemize}

Now the gravitational and effective potential are 
\begin{equation}
U=-\frac{1}{8}x^2+\frac{3\sqrt{3}}{4}(1-2\mu)xy+\frac{5}{8}y^2-\frac{1}{2}z^2+K(x,y,z),\label{hillgravpotential}
\end{equation}
\begin{equation}
\Omega=\frac{1}{2}(x^2+y^2)+U=\frac{3}{8}x^2+\frac{3\sqrt{3}}{4}(1-2\mu)xy+\frac{9}{8}y^2-\frac{1}{2}z^2+K(x,y,z),\label{hillefectivepotential}
\end{equation}
respectively. The equations of motion can be written as in the full problem
\begin{equation}
\label{hillfinalequations}
\ddot{x}-2\dot{y}=\Omega_{x},
\end{equation}
$$\ddot{y}+2\dot{x}=\Omega_{y},$$ $$\ddot{z}=\Omega_{z},$$
but $\Omega$ is given by the equation (\ref{hillefectivepotential}).

\section{The equilibrium points of the system.}
\subsection{Computation of the equilibrium points.}
In this section we prove that the system has 4 equilibrium points and we will be able to compute them explicitly in terms of the mass parameter $\mu$. So, in order to find the equilibrium points of the limit case, as usual, we need to find the critical points of the effective potential (\ref{hillefectivepotential}), an easy computation shows that $$\Omega_{z}=z(1+\frac{1}{r^3}),$$ the equation $\Omega_{z}=0$ implies that $z=0$ so the equilibrium points of the system are coplanar. Therefore, it is enough to study the critical points of the planar effective potential $$\Omega=\frac{3}{8}x^2+\frac{3\sqrt{3}}{4}(1-2\mu)xy+\frac{9}{8}y^2+K(x,y),$$ in matrix notation. 
\begin{equation}
\Omega=\frac{1}{2}z^{T}Mz+\frac{1}{\Vert z\Vert},
\end{equation}
where $z=(x,y)^{T}$ and $M$ is the matrix $$\left(\begin{array}{cc}\frac{3}{4}  & \frac{3\sqrt{3}}{4}(1-2\mu) \\\frac{3\sqrt{3}}{4}(1-2\mu)  & \frac{9}{4}\end{array}\right).$$
The above matrix has eigenvalues $$\lambda_{1}=\frac{3}{2}(1-d),$$ $$\lambda_{2}=\frac{3}{2}(1+d),$$ with respective eigenvectors $$v_{1}=\left(\frac{1+2d}{(2\mu-1)\sqrt{3+(\frac{1+2d}{1-2\mu})^2}},\frac{\sqrt{3}}{\sqrt{3+(\frac{1+2d}{1-2\mu})^2}}\right),$$ and $$v_{2}=\left(\frac{1-2d}{(2\mu-1)\sqrt{3+(\frac{1-2d}{1-2\mu})^2}},\frac{\sqrt{3}}{\sqrt{3+(\frac{1-2d}{1-2\mu})^2}}\right),$$ where $d=\sqrt{1-3\mu+3\mu^2}$. The eigenvectors have been chosen such that $\Vert v_{1}\Vert=\Vert v_{2}\Vert=1$. The equation to be solved is $\nabla\Omega=0$, or explicitly
\begin{equation}
Mz-\frac{z}{\Vert z\Vert^3}=0,\label{gradient}
\end{equation}
we can use the invertible matrix $C=col(v_{1},v_{2})$ to solve the above equation, if we consider the linear change of variables $z=Cz'$, substitute in the equation  (\ref{gradient}) and multiply by $C^{-1}$, we obtain $$C^{-1}MCz'-\frac{C^{-1}Cz'}{\Vert Cz\Vert^3}=0,$$ or equivalently 
\begin{equation}
Dz'-\frac{z'}{\Vert Cz'\Vert^3}=0,\label{gradientdiagonal}
\end{equation}
where $D$ is given by the diagonal matrix $$D=\left(\begin{array}{cc}\lambda_{1}  & 0 \\0 &\lambda_{2}\end{array}\right).$$ In terms of coordinates the equation (\ref{gradientdiagonal}) is equivalent to the system 
\begin{equation}
(\lambda_{1}-\frac{1}{\Vert Cz'\Vert^3})x'=0,
\end{equation}
\begin{equation}
(\lambda_{2}-\frac{1}{\Vert Cz'\Vert^3})y'=0,
\end{equation}
It is clear that in the above equations the case $x'=y'=0$ corresponds to a singularity and the case $x'\ne0$, $y'\ne0$ gives rise a contradiction, therefore when $y'=0$ we have $\Vert Cz'\Vert^3=\lambda^{-1}_{1}$  or equivalently $$\vert x'\vert=\frac{1}{\sqrt[3]{\lambda_{1}}\Vert v_{1}\Vert},$$
on the other hand, when $x'=0$ we have $\Vert Cz'\Vert^3=\lambda^{-1}_{2}$  or equivalently $$\vert y'\vert=\frac{1}{\sqrt[3]{\lambda_{2}}\Vert v_{2}\Vert},$$
but  $\Vert v_{1}\Vert=\Vert v_{2}\Vert=1$, therefore we obtain four equilibrium points given by 
\begin{eqnarray*}
L'_{1}=(0,\frac{1}{\sqrt[3]{\lambda_{2}}}), 
L'_{2}=(0,-\frac{1}{\sqrt[3]{\lambda_{2}}}),
L'_{3}=(\frac{1}{\sqrt[3]{\lambda_{1}}},0),
L'_{4}=(-\frac{1}{\sqrt[3]{\lambda_{1}}},0),
\end{eqnarray*}
or in the original coordinates we have $L_{i}=CL'^{T}_{i}$ for $i=1,2,3,4$. It is easy to see that 
\begin{eqnarray*}
L_{1}=\frac{1}{\sqrt[3]{\lambda_{2}}}v_{2}, 
L_{3}=\frac{1}{\sqrt[3]{\lambda_{1}}}v_{1},
\end{eqnarray*}
and $L_{2}=-L_{1}$, $L_{4}=-L_{3}$.

\subsection{Study of the stability of the equilibrium points.}
\begin{figure}[!hbp]
\centering
\includegraphics[width=3.5in]{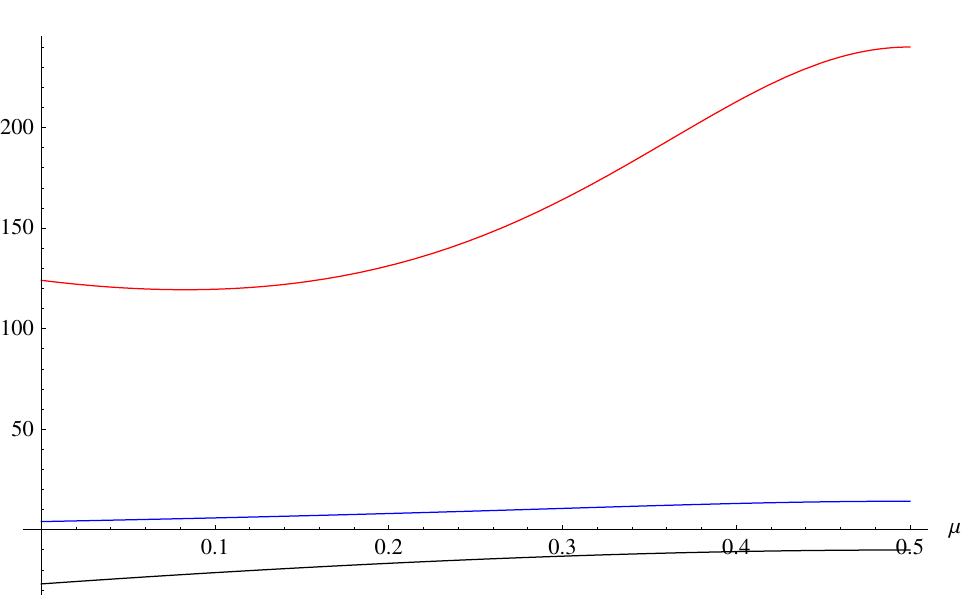}
\caption{The coefficients $A$ (in blue), $B$ (in black) and $D$ (in red) for the equilibrium point $L_{1}$ as functions of the mass parameter $\mu$.}\label{stabilityL1}
\end{figure}

\begin{figure}[!hbp]
\centering
\includegraphics[width=3.5in]{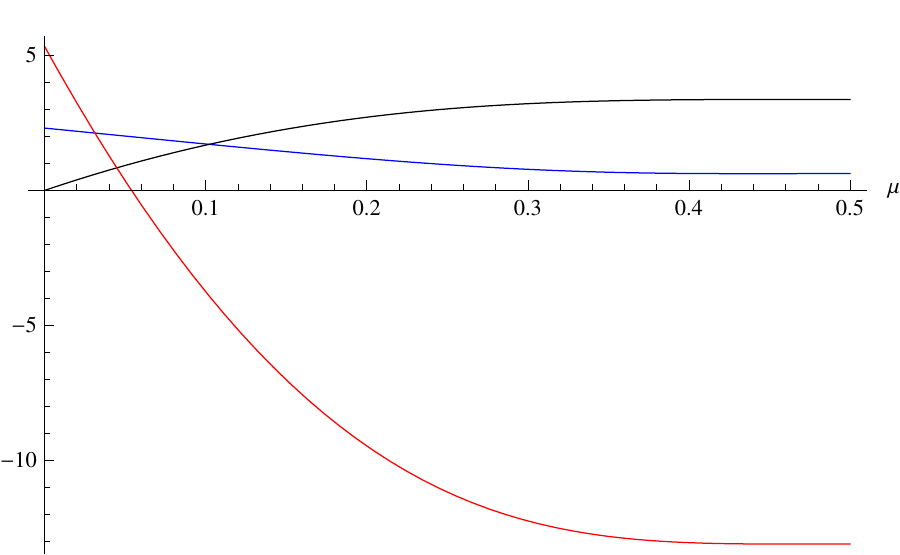}
\caption{The coefficients $A$ (in blue), $B$ (in black) and $D$ (in red) for the equilibrium point $L_{3}$ as functions of the mass parameter $\mu$.}\label{stabilityL3}
\end{figure}

In the previous subsection we obtained explicit expression of the four equilibrium points in terms of the parameter $\mu$ so, we can analyze the stability in the whole range $\mu\in[0,1/2]$, we will perform such analysis for the planar case $z=0$. We need to linearize the equations of motion, i.e., we need to study the linear system $\bf{\dot{\xi}}=\bf{A}\bf{\xi}$, where $\bf{\xi}=(x,y,\dot{x},\dot{y})^{T}$ and $A$ is the matrix 
\begin{equation}
\left(\begin{array}{cccc}0 & 0 & 1 & 0 \\0 & 0 & 0 & 1 \\\Omega_{xx} & \Omega_{xy} & 0 & 2 \\\Omega_{xy} & \Omega_{yy} & -2 & 0\end{array}\right)\label{matrixlinearization}
\end{equation}
where the partial derivatives $$\Omega_{xx}=\frac{3}{4}+\frac{3x^2}{(x^2+y^2)^{5/2}}-\frac{1}{(x^2+y^2)^{3/2}},$$ $$\Omega_{yy}=\frac{9}{4}+\frac{3y^2}{(x^2+y^2)^{5/2}}-\frac{1}{(x^2+y^2)^{3/2}},$$ and $$\Omega_{xy}=\frac{3\sqrt{3}}{4}(1-2\mu)+\frac{3xy}{(x^2+y^2)^{5/2}},$$ need to be evaluated in each $L_{i}$ for $i=1,2,3,4.$ However we must observe that because of the symmetry of the equilibrium points, we just need to study the equilibrium points $L_{1}$ and $L_{3}$. It is well known that the characteristic polynomial of the matrix (\ref{matrixlinearization}) is given by the expression 
\begin{equation}
p(\lambda)=\lambda^4+A\lambda^2+B,
\end{equation}
where $A=4-\Omega_{xx}-\Omega_{yy}$, $B=\Omega_{xx}\Omega_{yy}-\Omega_{xy}^2$. Therefore the four eigenvalues are $$\lambda_{1,2,3,4.}=\pm\frac{1}{\sqrt{2}}\sqrt{-A\pm\sqrt{D}},$$ with $D=A^2-4B$. A equilibrium point will be linearly stable if only if  $A$, $B$ and $D$ are non-negatives. Because of we know explicit expressions in term of the mass parameter $\mu$ for each equilibrium point, the referred coefficients of the characteristic polynomial are functions of $\mu$ so, using techniques of calculus we can study the behavior of such coefficients. In the figures \ref{stabilityL1} and \ref{stabilityL3} we can observe the behavior of  $A$, $B$ and $D$ as functions of the parameter $\mu$. 

\begin{prop} The coefficient $B$ in negative for $\mu\in[0,1/2]$ so the equilibrium point $L_{1}$ is unstable for this range of values of the mass parameter, in fact, the eigenvalues are given by $\pm\lambda$ and $\pm\textit{i}\omega$ with $\lambda>0$ and  $\omega>0$.
\end{prop}

On the other hand, we can observe that the discriminant $D$ changes form positive to negative for the equilibrium point $L_{3}$, therefore
\begin{prop} There exists a value $\mu_{0}$ such that $D=0$, as a consequence, this equilibrium point has the following properties: for $\mu\in[0,\mu_{0})$ the eigenvalues are $\pm\textit{i}\omega_{1}$ and $\pm\textit{i}\omega_{2}$, for $\mu=\mu_{0}$ we have a pair of the eigenvalues $\pm\textit{i}\omega$ of multiplicity 2, finally when $\mu\in(\mu_{0},1/2]$ the eigenvalues are $\pm\alpha\pm\textit{i}\omega$ with $\alpha>0$ and  $\omega>0$.
\end{prop}

\begin{figure}
  \centering
\begin{tabular}{ccc}
  \includegraphics[width=1.7in]{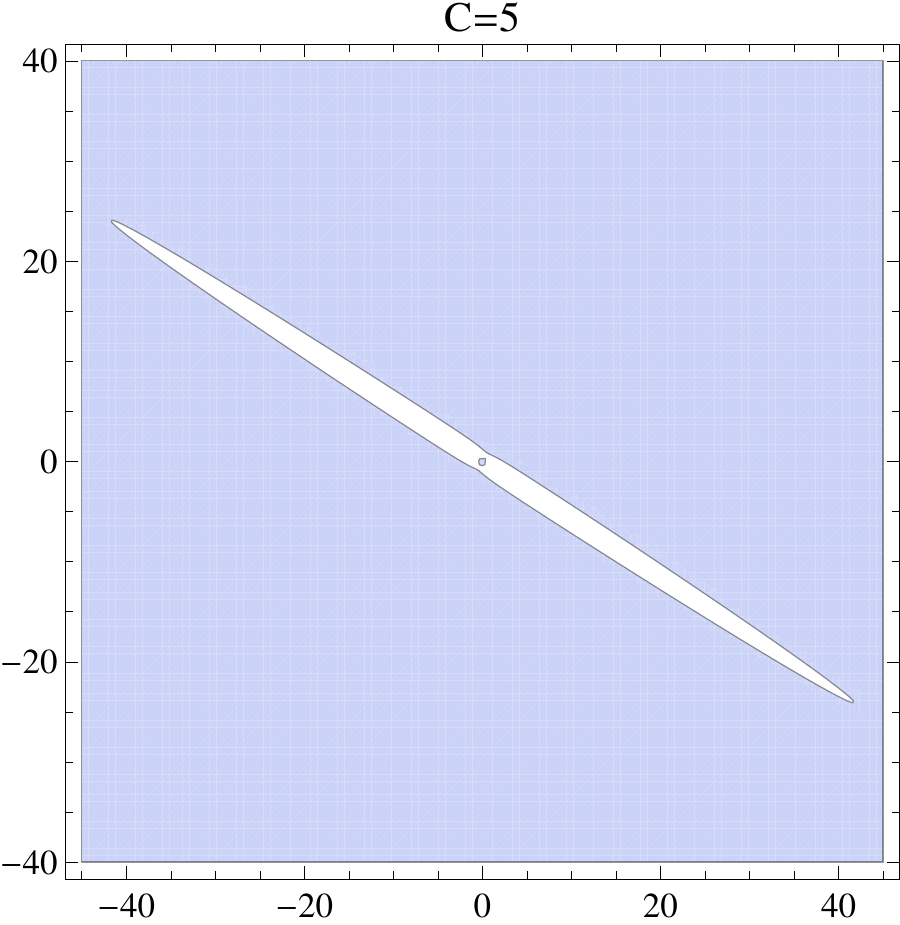}& \includegraphics[width=1.7in]{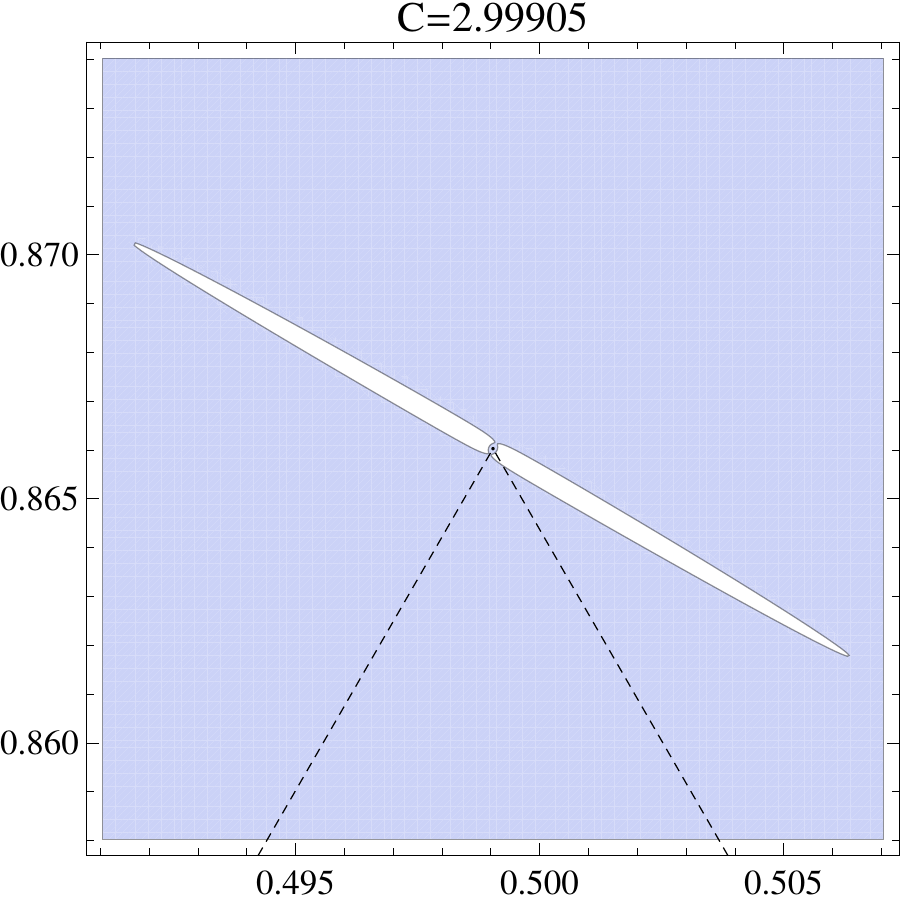}& \includegraphics[width=1.7in]{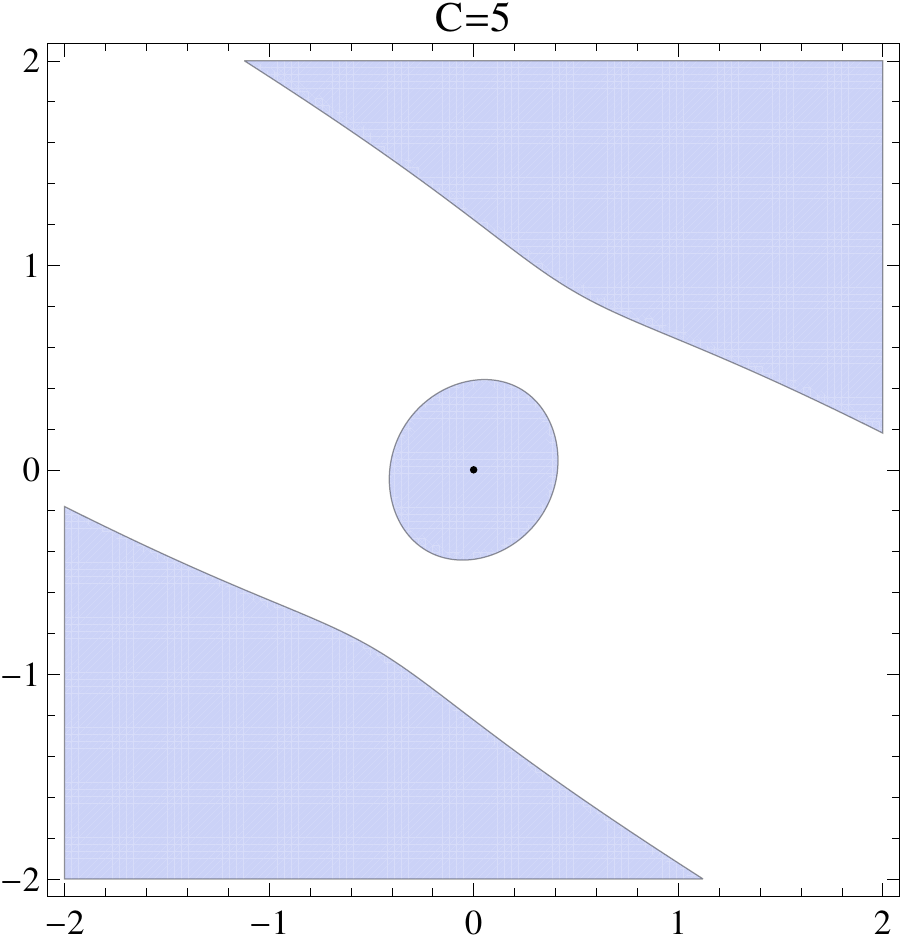} \\
  \includegraphics[width=1.7in]{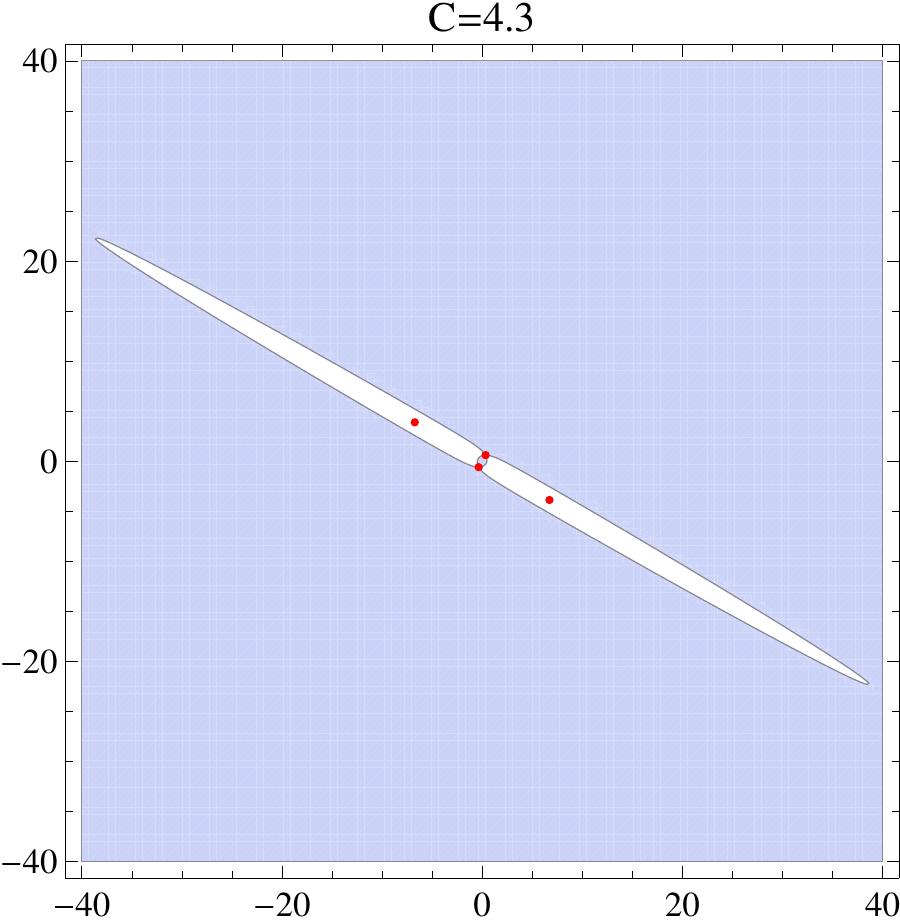}& \includegraphics[width=1.7in]{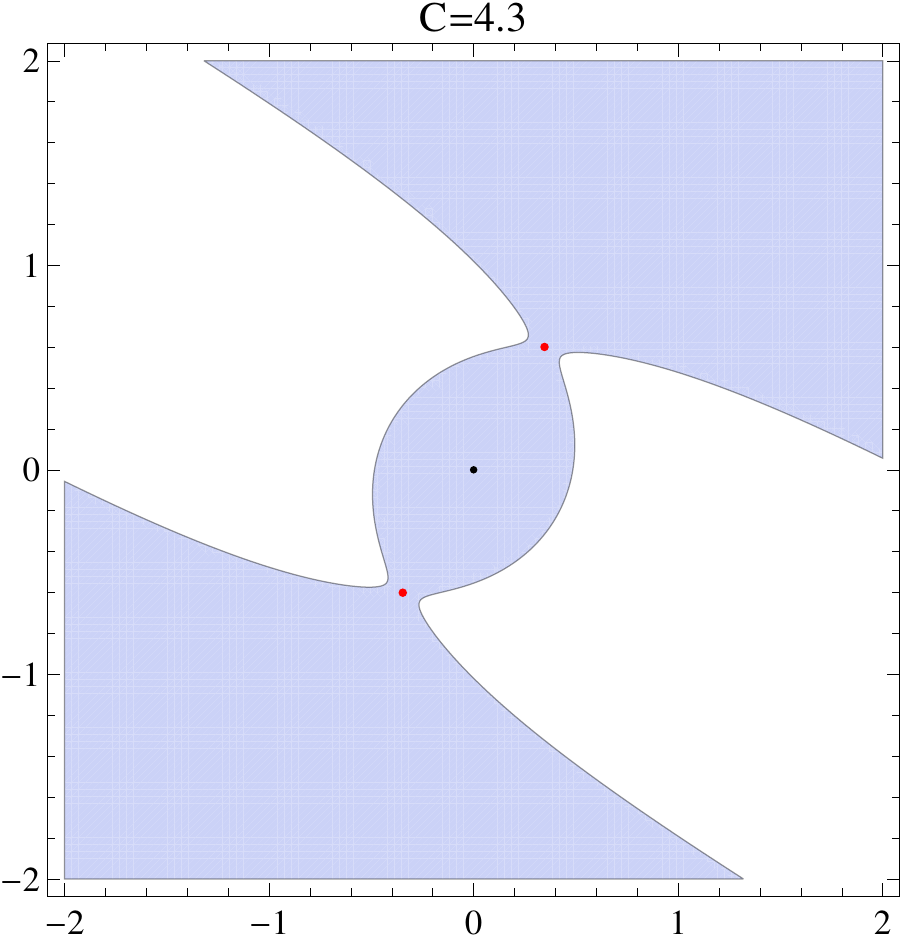}& \includegraphics[width=1.7in]{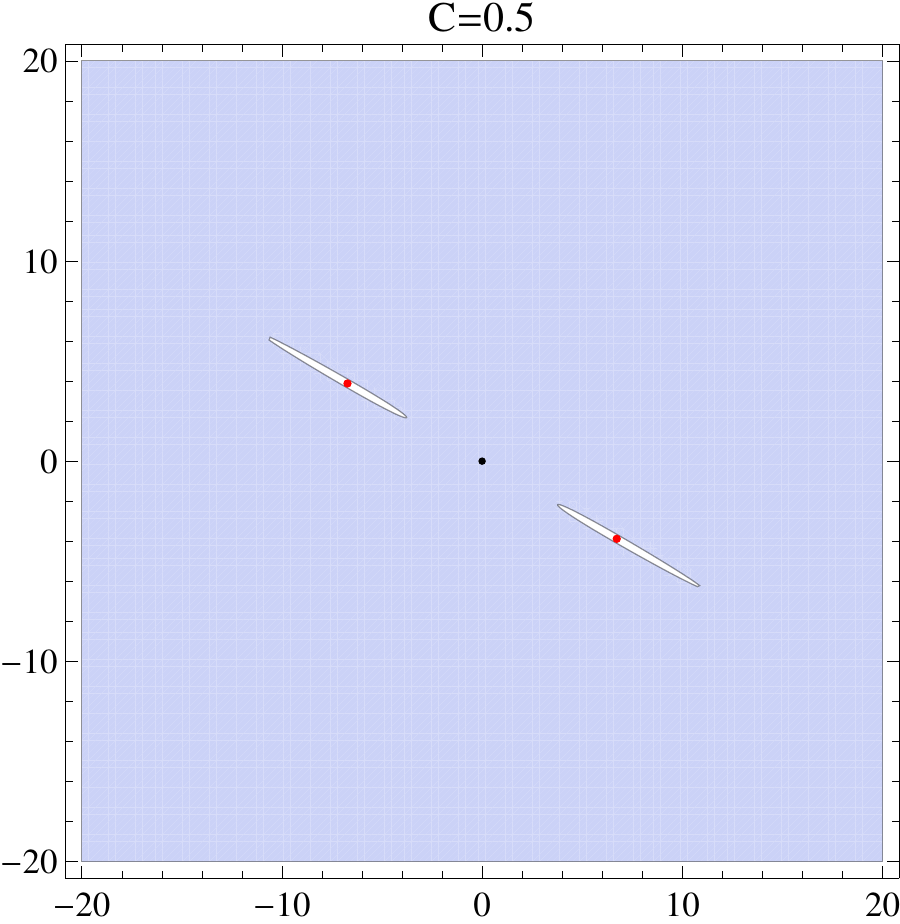}
\end{tabular}
 \caption{Hill's regions (blue areas) for for $\mu= 0.00095$. First row, left to right, Hill's region for a the limit problem, Hill's region for a the full R4BP when $m_{3}=7.03\times10^{-12}$ (mass of 624 Hektor), magnification of the first figure for the limit case. Second row, left to right, position of the equilibrium points (red dots) for the limit case.}\label{limithillregions}
\end{figure}

In the figure \ref{limithillregions} we show the so called Hill's regions for the planar case and for $\mu= 0.00095$ that corresponds to mass ratio of the Sun-Jupiter system, in the first two figures of the first row we show the Hill's regions for the limit problem and for the full R4BP when $m_{3}=7.03\times10^{-12}$, the mass ratio of the asteroid 624 Hektor, the lines in the second figure are imaginary lines that connect $m_{3}$ with the remaining masses. We have marked the position of the fixed mass with a black dot and the positions of the four equilibrium points with red dots.

\newpage

\begin{thebibliography}{99}
\bibitem{PapaIII} Baltagiannis, A.N., Papadakis, K.E.; Periodic solutions in the Sun-Jupiter-Trojan Asteroid-Spacecraft system. Planetary and Space Science. \textbf{75} ,148--157 (2013).
\bibitem{Cecc} Ceccaroni M., Biggs J.; Extension of low-thrust propulsion to the autonomous coplanar circular restricted four body problem with application to future Trojan Asteroid missions. In: 61st Int. Astro. Congress IAC 2010 Prague, Czech Republic (2010).
\bibitem{Hill} Hill G.W., Researches in the Lunar Theory.American Journal of Mathematics. \textbf{1}. No. \textbf{1} 5-26 (1878).
\bibitem{Marchis}Marchis, F., et al: The Puzzling Mutual Orbit of the Binary Trojan Asteroid (624) Hektor.  Astropysical Letters, ApJ \textbf{783}, L37. (2014).
\bibitem{MeyerHDS} Meyer K.; Introduction to Hamiltonian Dynamical Systems and the N-body problem. Springer Verlag. (2009)
\end{thebibliography}
\end{document}